# An elementary non-recursive expression for the partition function P(n).

Godofredo Iommi Amunátegui

Instituto de Física

Pontificia Universidad Católica de Valparaíso

giommi@ucv.cl

Abstract

Consideration of a classification of the number of partitions of a natural number according to the members of sub-partitions differing from unity leads to a non-recursive formula for the number of irreducible representations of the symmetric group $S_n$.

This article was published, long ago, under the title "A non-recursive expression for the number of irreducible representations of the Symmetric Group $S_n$, Physica 114A, 1982, 361-364, North-Holland Publishing Co. The Introduction has been, somewhat, improved, however, the handmade result remains unproved.

1. **Introduction**

As is well known, there are as many inequivalent irreducible representations of the symmetric group $S_n$ as there are partitions of $n$. Indeed the concept of Young diagrams, which is most useful for discussing the basis functions of the unitary irreducible representations of $SU_n$, is introduced through the consideration of $S_n$. Or, in Weyl's words[1], the substratum of a representation of the unitary group consists of the linear manifold of all tensors of order $n$ whose symmetry properties are expressed by linear relations between them and the tensors obtained from them by the $n!$ permutations. To each Young diagram corresponds a partition of $n$. Hardy and Ramanujan[2] obtained an asymptotic formula for the number of partitions of $n$, $p(n)$, of noteworthy accuracy. Hans Rademacher, following their steps, converted such an approximation into an identity[3]. More recently Conway and Guy have complained about the fact that "There's no simple exact formula for $p(n)$"[4], and Bruce Sagan, in his wonderful book on the Symmetric Group[5], apropos the generating function for partitions due to Euler, remarks: "Despite the simplicity of this generating formula, there is no known closed formula for $p(n)$ itself".

Tables are available which list $p(n)$ for extensive ranges of $n$. These results are obtained from the recursion formula[6] $p(n) = \sum (-1)^{j+1} p\left(n - \tfrac{1}{2}(3j^2 + j)\right)$, where the sum runs over all the positive integers such that the arguments of the partition function are non-negative. It would therefore be interesting to have a direct, non-recursive algorithm for $p(n)$ and in this note we present a fairly elementary procedure which leads to such an expression.

2. **Derivation of a simple expression for $p(n)$.**

Let $p(n)$ be the number of partitions of $n$, i.e., the number of ways of representing $n$ as a sum of positive integral parts. The order in which the parts are arranged is deemed to be irrelevant. In table I, the partitions are classified according to the numbers of their sub-partitions which are different from unity. Thus, the partition $2^2 1^2$ of $n = 6$ has two such partitions while the partition $2^3$ would have three.

Table I: Classification of the partitions of *n* according to the numbers of sub-partitions differing from unity

| p(n) | Numbers of sub-partitions different from unity | | | | | |
|---|---|---|---|---|---|---|
|  | 0 | 1 | 2 | 3 | 4 | 5 |
| p(1) | 1 | 0 | 0 | 0 | 0 | 0 |
| p(2) | 1 | 1 | 0 | 0 | 0 | 0 |
| p(3) | 1 | 2 | 0 | 0 | 0 | 0 |
| p(4) | 1 | 3 | 1 | 0 | 0 | 0 |
| p(5) | 1 | 4 | 2 | 0 | 0 | 0 |
| p(6) | 1 | 5 | 4 | 1 | 0 | 0 |
| p(7) | 1 | 6 | 6 | 2 | 0 | 0 |
| p(8) | 1 | 7 | 9 | 4 | 1 | 0 |
| p(9) | 1 | 8 | 12 | 7 | 2 | 0 |
| p(10) | 1 | 9 | 16 | 11 | 4 | 1 |
| p(11) | 1 | 10 | 20 | 16 | 7 | 2 |

Addition of the terms in the first two columns of table 1 gives *n*. Let $A^0_n$ denote this sum, $A^1_n$ the third column terms, $A^2_n$ the fourth column terms and so on. We may therefore rewrite table I in the form of table II.

Table II: Table of the functions $A^j_n$ for values of *n* in the range $1 \leq n \leq 11$

| n | $A^0_n$ | $A^1_n$ | $A^2_n$ | $A^3_n$ | $A^4_n$ |
|---|---|---|---|---|---|
| 1 | 1 | 0 | 0 | 0 | 0 |
| 2 | 2 | 0 | 0 | 0 | 0 |
| 3 | 3 | 0 | 0 | 0 | 0 |
| 4 | 4 | 1 | 0 | 0 | 0 |
| 5 | 5 | 2 | 0 | 0 | 0 |
| 6 | 6 | 4 | 1 | 0 | 0 |
| 7 | 7 | 6 | 2 | 0 | 0 |
| 8 | 8 | 9 | 4 | 1 | 0 |
| 9 | 9 | 12 | 7 | 2 | 0 |
| 10 | 10 | 16 | 11 | 4 | 1 |
| 11 | 11 | 20 | 16 | 7 | 2 |

By examining table II we can see that certain relationships between the $A^j_n$ are valid for all values of $n$, for example $A^2_{n+5} = A^1_{n+3} + A^1_n$. A similar analysis can be carried out for the entries of any column. In general, therefore, the $i$th column terms can be written as a sum of $A^1_n$ terms. This procedure enables us to formulate the following non-recursive expression for $p(n)$:

$$p(n) = \sum_{\beta=0}^{r} A^\beta_n,$$

where $r = \tfrac{1}{2}(n-2)$, $n$ even; $r = \tfrac{1}{2}(n-3)$, $n$ odd.

When $\beta \geq 2$, the auxiliary functions $A^\beta_n$ are given in general by the formula

$$A^\beta_n = \sum_{\alpha_1=0}^{s-2} \sum_{\alpha_2=0}^{s-3} \cdots \sum_{\alpha_{\beta-1}=0}^{s-\beta} A^1_{n+2-2\beta-\gamma},$$

where $\gamma = \sum_{\varepsilon=0}^{\beta-2}(\varepsilon+3)\alpha_{\varepsilon+1}$, $s = \tfrac{1}{3}n \pmod 3$ and $\alpha_{j-1} = 0$ when $s \leq j$.

When $\beta = 1$, the appropriate formula is

$$A^1_n = \begin{cases} r^2 & n \text{ even}, \\ r(r+1), & n \text{ odd}, \end{cases}$$

while when $\beta = 0$, $A^0_n$ is, by definition, equal to $n$.

As an example, we shall calculate $p(22)$; when $n = 22$, $r = 10$ and $s = 7$.

Hence

$$p(22) = \sum_{\beta=0}^{10} A^\beta_{22} = A^0_{22} + A^1_{22} + A^2_{22} + A^3_{22} + A^4_{22} + A^5_{22} + A^6_{22} + A^7_{22} + A^8_{22} + A^9_{22} + A^{10}_{22}.$$

Now $A^0_{22} = 22$ and $A^1_{22} = 10^2 = 100$: the remaining terms are obtained from the formula as

$$A^2_{22} = \sum_{\alpha_1=0}^{5} A^1_{20-3\alpha_1} = A^1_{20} + A^1_{17} + A^1_{14} + A^1_{11} + A^1_8 + A^1_5 = 204,$$

$$A^3_{22} = \sum_{\alpha_1=0}^{5}\sum_{\alpha_2=0}^{4} A^1_{18-3\alpha_1-4\alpha_2} = A^1_{18} + A^1_{15} + A^1_{12} + A^1_9 + A^1_6 + A^1_{14} + A^1_{11} + A^1_8 + A^1_5 + A^1_{10} + A^1_7 + A^1_4 + A^1_6 = 241,$$

$$A^4_{22} = \sum_{\alpha_1=0}^{5}\sum_{\alpha_2=0}^{4}\sum_{\alpha_3=0}^{3} A^1_{16-3\alpha_1-4\alpha_2-5\alpha_3} = 197.$$

By an analogous method, we have

$A^5_{22} = 125$, $A^6_{22} = 66$, $A^7_{22} = 30$, $A^8_{22} = 12$, $A^9_{22} = 4$ and $A^{10}_{22} = 1$,

So that, finally, we obtain

$$p(22) = 22 + 100 + 204 + 241 + 197 + 125 + 66 + 30 + 12 + 4 + 1 = 1002.$$

---

[1] H. Weyl, The Theory of Groups and Quantum Mechanics, Dover Pub., p. 281.
[2] G.H. Hardy and S. Ramanujan, Compt. Rend. Acad. Sci. (Paris) 164 (1917) p.35.
[3] H. Rademacher, On the partition function $p(n)$. Proc. of the London Math. Soc., 1938, Vol. 2, Nº 1, pp. 241-254.
[4] J.H. Conway, R.K. Guy, The Book of Numbers, Springer 1996, p. 95.
[5] B. Sagan, The Symmetric Group: Representations, Combinatorial Algorithms, and Symmetric functions. Springer, second edition, 2010, p. 144.
[6] This formula may, also, be traced to Euler.